\documentclass{article}
\usepackage{amsmath}
\font\tengoth=eufm10
\font\sevengoth=eufm7
\font\fivegoth=eufm5
\newfam\gothfam
\textfont\gothfam=\tengoth \scriptfont\gothfam=\sevengoth
\scriptscriptfont\gothfam=\fivegoth

\def\blacksquare{\hbox to .60em{\vrule width .60em height .60em}}
  \font\bb=msbm10 
\catcode`\@=12  
\def\é{\'e}
\def\{\`e}
\def\?{\`a}
\def\{\`u}
\def\{\c c}
\def\hb {\hfil \break}
\def\n {\vskip 0.2cm \noindent }
\def\scirc{\,{\raise 0.8pt\hbox{$\scriptstyle\circ$}}\,}
\def\ins{\,{\raise 0.2cm \hbox{ $\scriptstyle \circ$}}\,}

\def  \é{\'e}
\def\è{\`e}
\def\à{\`a}
\def\ù{\`u}
\def\ç{\c c$\!\!\!$}

\date{}
 
\begin{document}
  
 \centerline {\bf  Courbure  des  tissus planaires d\éfinis implicitement } 
  \centerline {\bf  par une \équation diff\érentielle polynomiale en $y'$.}
  \centerline {\bf Programmation}
%\centerline {    d'apr\ès A. H\énaut}
 
 \bigskip
 
 \rightline {\bf   J.P. Dufour et D. Lehmann}  
  \bigskip

%\rightline {Etat du 30/04/2015}
 \bigskip
 
{ \bf Abstract : }

{\it The  aim of this paper is mainly, after some theoretical explanations,  to provide a program on Maple for computing, whatever be $d$, the 
curvature of the planar $d$-web implicitely defined by a  differential equation $F(x,y,y')=0$, $F$ being polynomial of degree $d$ with respect to $y'$.

 Moreover, we prove in the appendix a "concentration theorem" for  any calibrated ordinary $d$-web of codimension one in a $n$-dimensional manifold $($in particular for any planar web$)$.  Its curvature matrix,   relatively  to an "adapted" trivialization, is  concentrated on the $(n-2+k_0)!/(n-2)!k_0!$ last lines $($the last line if $n=2)$, $k_0$ denoting the integer $\geq 2$ such that $d=(n-1+k_0!/(n-1)!k_0!$.

}
 
\section{Introduction }

A. H\énaut  a beaucoup \étudi\é les tissus planaires d\éfinis par une  \équation diff\érentielle $F(x,y,y')=0$, o\ù 
$$F(x,y,p)=\sum_{i=0}^d  a_i(x,y)\ p^{d-i}$$
d\ésigne un polyn\^ome de degr\é $d$ en $p$, \à coefficients $a_i$ dans l'anneau   ${\cal O}_U$  des fonctions holomorphes sur  un ouvert $U$ de {\bb C}$^2$.

 On suppose
 
- que, pour tout $(x,y)$ dans $U$, $F(x,y,p)$  et $F'_p(x,y,p)$  ne sont jamais simultan\ément nuls, c'est-\à-dire que les racines $p_i(x,y)$ $(i=1,\cdots, d)$ de $F$,  consid\ér\é comme polyn\^ome en $p$, sont toutes distinctes, 

- et, pour simplifier les  calculs,  que  le coefficient $a_0 $ de $p^d$ est identiquement \égal \à 1, soit
 $$F(x,y,p)\equiv\prod_{i=1}^d \bigl(p-p_i(x,y)\bigr).$$  Notons ${\cal F}_i$ le feuilletage d\éfini par l'\équation diff\érentielle $y'=p_i(x,y)$ ;    en superposant les $d$ feuilletages ${\cal F}_i$, on obtient un $d$-tissu holomorphe d\écomposable sur $U$. 

Dans [H] en particulier, A. H\énaut a d\éfini, pour un tel tissu, une certaine connexion dont la courbure  g\én\éralise \à tout $d$ la courbure de Blaschke-Dubourdieu ([B]) du cas $d=3$, et qui est l'obstruction \à ce que le rang du tissu (c'est-\à-dire la dimension de l'espace vectoriel de ses relations ab\éliennes) soit maximum, \égal \à $\frac{1}{2}(d-1)(d-2)$.

%\footnote{Dans le cas $d=3$ au moins,  les r\ésultats donn\és par le programme Ripoll sont erron\és :   le 3-tissu d\éfini  par $a_0=-a_2=1$ et $a_1=-a_3=G.M'_x(x,y)/(1+G.M'_y(x,y)) $  \à l'aide  d'une fonction $M(x,y)$ de deux variables et d'un param\ètre scalaire $G$, est une d\éformation $\bigl(x+y,\ x-y,\ y+G.M(x,y)\bigr)$   du 3-tissu hexagonal $(x+y,x-y,y)$ par le param\ètre $G$ : sa courbure, qui doit \^etre nulle pour   $G=0$, est de la forme    $(M'''_{x^3}-M'''_{x,y^2})G+O(G^2)$, alors  qu'avec  le programme Ripoll, on trouve  $(M'''_{x,y^2}-(3/2)M'''_{x^3})G+O(G^2)$, les termes en $G^2$ \étant encore plus disparates.} 

Dans sa th\èse, \à 
la m\^eme \'epoque,  L. Pirio \cite{Pi} a en particulier modernis\é des travaux anciens  de A. Pantazi \cite{Pa}, é \également en termes de connexion.

Des programmes  permettant de calculer cette courbure ont  \ét\é r\édig\és sur Maple,  pour $d=3,\ 4$ ou $5$,  dans le cas "explicite"  par Pirio ( \cite{Pi}), et dans le cas "implicite"  par O. Ripoll (\cite{R1}).

Dans [CL], nous avons g\én\éralis\é la d\éfinition de cette courbure aux $d$-tissus de codimension un  sur une  vari\ét\é holomorphe de dimension $n$ arbitraire $(n\geq 2)$, pourvu que ceux-ci soient ``ordinaires'' (condition  toujours v\érifi\ée pour $n=2$, et g\én\ériquement v\érifi\ée localement pour $n\geq 3$), et ``calibr\és''

\pagebreak \n (ce qui veut dire qu'il existe un entier $k_0$ ($k_0\geq 2)$, tel que $d=\frac{(n-1+k_0)!}{(n-1)!k_0!}$, cette condition \étant  toujours v\érifi\ée pour $n=2$ avec $k_0=d-1$). 

Dans [DL], nous avons r\édig\é un programme sur Maple, permettant de calculer  cette courbure quels que soient $n$ et $k_0$, en un temps raisonnable pourvu que $n$ et $k_0$ ne soient pas ``trop grands''. 
Toutefois, pour appliquer ce programme, il fallait disposer d'une int\égrale premi\ère pour chaque feuilletage ${\cal F}_i$, 
 et par  cons\équent, lorsque $n=2$,  supposer int\égr\ée chaque \équation diff\érentielle $y'=p_i(x,y)$. En fait, nous disposons aussi d'un autre programme (non publi\é) qui s'applique au cas o\ù chaque feuilletage local du tissu est d\éfini par une forme int\égrable non n\éc\éssairement ferm\ée, ce qui -pour $n=2$, revient  au cas  "explicite" o\ù l'on conna\^it au moins   la d\écomposition  $F(x,y,p)\equiv \prod_{i=1}^d\bigl(p-p_i(x,y)\bigr)$ . Dans le cas "implicite" 
o\ù l'on ne sait pas r\éaliser cette d\écomposition, il nous faut g\én\éraliser \à tout $d$ une  m\éthode de programmation du type de celle propos\ée par Ripoll dans les cas $d=3,4,5$. C'est l'objet de ce travail\footnote{Dans le cas explicite o\ù l'on a le choix, le programme explicite est en g\én\éral plus performant.}. 

 Nous r\ésumons ci-dessous les quelques points de la th\éorie dont nous avons besoin pour expliquer la programmation, et renvoyons \à  [Pi], [H], [CL] ou [DL] pour plus de d\étails.

\section{Rappels sur la d\éfinition des relations ab\éliennes} 

Il revient au m\^eme de se donner la diff\érentielle d'une int\égrale premi\ère du feuilletage ${\cal F}_i$ ou une fonction holomorphe $f_i:(x,y)\mapsto f_i(x,y)$ telle que la 1-forme $f_i(dy-p_idx)$ soit ferm\ée. Une relation 
ab\élienne du $d$-tissu obtenu par superposition des ${\cal F}_i$ est alors d\éfinie par une famille de
fonctions $f_i$ ($i=1,\cdots,d$), telle que

- la somme $\sum_{i=1}^d f_i(dy-p_idx)$ soit nulle (condition dite ``de trace nulle''), 

- chaque forme $f_i(dy-p_idx)$ est ferm\ée, soit :

$(I)$\hskip 1cm $\sum_i f_i\equiv 0$, 

$(II)$\hskip .8cm $\sum_i p_i  f_i\equiv 0$, 

$(III_i)$\hskip .6cm   $(f_i)'_x +(p_if_i)'_y\equiv 0$ pour tout $i$.

Notons $\tilde U$ l'ouvert de la vari\ét\é de contact au dessus de $U$, admissible pour les coordonn\ées $(x,y,p)$, et $S$ la surface  d'\équation $F(x,y,p)=0$ dans $\tilde U$. Soit $ {\cal F} $ le feuilletage sur  $S$ d\éfini par la restriction  $(dy-p\ dx)_{_S} $ \à cette surface
 de la forme de contact $dy-p\ dx$ sur $\tilde U$ : la surface $S$ est un rev\^etement (trivial) \à $d$ feuillets de $U$.   Notant 
  $\pi_i:U_i\buildrel\cong \over\rightarrow U$ la projection sur $U$ du feuillet   $U_i$  d'\équation $p=p_i(x,y)$, la restriction $\tilde {\cal F}_i $ de $\tilde {\cal F} $ \à $U_i$ se projette sur  ${\cal F}_i $ par $\pi_i$.

 Soit $(f_i)_i$ une famille de fonctions sur $U$. Il revient au m\^eme de se donner une fonction $  f$    sur $S$ (celle dont la restriction \à $U_i$ est \égale \à $f_i\scirc\pi_i$). Et il revient aussi au m\^eme de dire que la forme $  f\ (dy-p\ dx)_{_S}$ sur $S$ est ferm\ée ou que chaque forme $f_i(dy-p_i\ dx)$ l'est. Mais exprimer cette condition, ainsi que la nullit\é 
 de la trace  pose probl\ème,  quand on ne sait pas identifier chaque feuillet $p=p_i(x,y)$. 
 
 \subsection{La m\éthode de H\énaut :}
 
Elle   consiste \à
 prolonger    la fonction $f:S\to ${\bb C} en une fraction rationnelle sur $\tilde U$ de la forme $$\tilde f(x,y,p)=\frac{r(x,y,p)}{F'_p(x,y,p)}$$ (donc sans p\^ole sur $S$), $r$ d\ésignant le  polyn\^ome d'interpolation en $p$ des expressions $f_i(x,y).F'_p\bigl(x,y,p_i(x,y)\bigr)$. Son  
  degr\é est a priori  au plus $ d-1$, mais nous allons voir qu'il est  m\^eme au plus $ d-3$ d\ès lors que sont v\érifi\ées les relations $(I)$ et $(II)$ :$$ r(x,y,p)=\sum_{j=0}^{d-3} r_j(x,y) p^j.$$ On prend alors   comme fonctions inconnues d\éfinissant les relations ab\éliennes  les $d-2$ coefficients  $r_j$ de  $r$, ($j=0,\cdots,d-3$), au lieu des $d$ fonctions $f_i$ ($i=1,\cdots,d$) reli\ées par les relations $(I)$ et $(II)$.
 
 Notons ${\cal O}_k[p]$ le ${\cal O}_U$-module des polyn\^omes $P\in  {\cal O}_U [p]$, de degr\é $\leq k$, et soit $$L_i(x,y,p)=\prod _{j,j\neq i}\bigl(p-p_j(x,y)\bigr)$$ dans ${\cal O}_{d-1}[p]$ : puisque $F(x,y,p)=\prod _{ i}\bigl(p-p_i(x,y)\bigr)$, $F'_p$ est 
 \égal \à $\sum_i L_i$. Posons  $$r=\sum_i f_i\ L_i.$$
 Puisque $\frac{r}{F'_p}$ est \égal \à $f_i$ sur le feuillet $U_i$,  $\frac{r}{F'_p}$ prolonge $f$ \à un voisinage de $S$ dans $\tilde U$.  
 
 A priori, le degr\é de $r$ est \égal \à  $d-1$. Mais le coefficient $r_{d-1}$ est \égal \à $\sum_i f_i$, tandis que $r_{d-2}=\sum_i p_i f_i+ a_1(\sum_i f_i).$ Les relations $(I)$ et $(II)$ deviennent donc $r_{d-1}=0$, et $r_{d-2}=0$, c'est-\à-dire : $r\in {\cal O}_{d-3}[p]$.
 
% Plus g\én\éralement, on obtient la formule : 
 
% \n  {\bf Lemme 1 : }
 
% $$r_j=  \sum_{k=0}^j a_k\sum_{i =1}^d f_i(p_i)^{j-k}

 \n {\bf Remarque :} Au lieu de prolonger la fonction  $f$ par  la fraction rationnelle $\frac{r}{F'_p}\ \Bigl(= \frac{\sum_i f_i\ L_i}{\sum_i L_i}\Bigr)$, on aurait pu la prolonger par  le polyn\^ome d'interpolation $\sum_i \frac{f_i}{ L_i(p_i)}L_i=\sum_{j=0}^{d-1} t_j(x,y)\ p^j$ des $f_i$, et prendre comme fonctions inconnues les $d$ coefficients $t_j$ reli\és par les deux \équations $(I)$ et $(II)$ (qui s'expriment facilement \à l'aide des seuls coefficients $a_i$ de $F$). Mais cette m\éthode est plus compliqu\ée du point de vue de la programmation : si l'on r\ésume les \équations $(I)$ et $(II)$ sous la forme 
 $$(**)\hskip 1cm <TT(a), t>\equiv 0,$$ ($TT(a)$ d\ésignant une matrice $2\times d$ construite\footnote{Pour $d=3$ par exemple, $TT(a)=\begin{pmatrix}(a_1)^2-2a_2&-a_1&3\\-(a_1)^3+3a_1a_2-3a_3&(a_1)^2-2a_2&-a_1 \end{pmatrix}$.} \à partir des coefficients $a_i$), il faut prendre les d\ériv\ées partielles successives de cette identit\é jusqu'\à l'ordre $d-2$, prendre le noyau du syst\ème obtenu, d\éfinir une trivialisation de ce noyau, etc... La m\éthode H\énaut est plus simple, puisque le nombre des fonctions inconnues $r_j$ est d\éj\à r\éduit \à $d-2$. 
 
 \n Il  reste maintenant  \à exprimer que la restriction \à $S$ de la 1-forme $\frac{r(x,y,p)}{F'_p(x,y,p)}(dy-p\ dx)$ est ferm\ée.

 \n {\bf Lemme 1 :} {\it Sur la surface  $S$, est v\érifi\ée l'\égalit\é $$d\Bigl( \frac{r(x,y,p)}{F'_p(x,y,p)}(dy-p\ dx)\Bigr)=\Biggl((r'_x+p\ r'_y)-\biggl(\frac{F'_x+p\ F'_y}{F'_p}\biggr)'_p.\ r-\frac{F'_x+p\ F'_y}{F'_p}\ r'_p\Biggr)\frac{dx\wedge dy }{  F'_p }.$$ 
}
 \n {\it D\émonstration :}  Posons $K:=F'_x+p\ F'_y$, d'o\ù $F_{x p}''+p\ F_{y p}''=K'_p-F'_y.$
 
  \n Sur $\tilde{U}$,  $d\Bigl(\frac{r(x,y,p)}{F'_p(x,y,p)}(dy-p\ dx)\Bigr)= d\Bigl(\frac{r(x,y,p)}{F'_p(x,y,p)}\Bigr)\wedge (dy-p\ dx)+\frac{r(x,y,p)}{F'_p(x,y,p)}(-dp\wedge dx),$ soit :\hb 
$$\Biggl[\biggl(\frac{r(x,y,p)}{F'_p(x,y,p)}\biggr)'_x+p\biggl(\frac{r(x,y,p)}{F'_p(x,y,p)}\biggr)'_y\Biggr]\ dx\wedge dy+\biggl(\frac{r(x,y,p)}{F'_p(x,y,p)}\biggr)'_p dp\wedge (dy-p\ dx)-\frac{r(x,y,p)}{F'_p(x,y,p)}(dp\wedge dx). $$
 L'expression $\Biggl[\biggl(\frac{r(x,y,p)}{F'_p(x,y,p)}\biggr)'_x+p\biggl(\frac{r(x,y,p)}{F'_p(x,y,p)}\biggr)'_y\Biggr]$
est \égale \à  $\frac{(r'_x+p\ r'_y)F'_p-r(K'_p-F'_y)}{(F'_p)^2}$. 

\n D'autre part, $dF\equiv 0$ sur la surface $S$ d'\équation $F(x,y,p)=0$, d'o\ù   $dp=\frac{-1}{F'_p} (F'_x\ dx +\ F'_y\ dy)$, 
$dp\wedge dx=\frac{F'_y}{F'_p} \ dx\wedge dy$ et $dp\wedge dy=-\frac{F'_x}{F'_p} \ dx\wedge dy$. Ainsi, sur $S$, 
$d\Bigl(\frac{r(x,y,p)}{F'_p(x,y,p)}(dy-p\ dx)\Bigr)=A \ dx\wedge dy,$  
$$\begin{matrix}\hbox{avec \hskip 1cm }A&=& \frac{(r'_x+p\ r'_y)F'_p-r.K'_p}{(F'_p)^2}-\biggl(\frac{r(x,y,p)}{F'_p(x,y,p)}\biggr)'_p.\frac{K}{F'_p} ,\\ 
&=&\frac{1}{(F'_p)^3}\Bigl((r'_x+p\ r'_y)(F'_p)^2+r (K.F_{p^2}''-K'_p.F'_p)-r'_p .K.F'_p)\Bigr),\\
&=&\frac{1}{(F'_p)^3}\Bigl((r'_x+p\ r'_y)(F'_p)^2+H.r-L.r'_p)\Bigr),
\end{matrix}$$
o\ù l'on a pos\é :$$H:=-(F'_p)^2.(K/F'_p)'_p\ (=K.F''_{p^2}-K'_p.F'_p) ,\hskip 1cm\hbox{  et }L:=K.F'_p.$$ 
  Il r\ésulte de ce lemme que la 1-forme $\frac{r(x,y,p)}{F'_p(x,y,p)}(dy-p\ dx)$ sera ferm\ée sur $S$ ssi :\hb {\it il existe $c\in {\cal O}_{2d-4}\ [p]$ tel que }$$(*_0)\hskip 1cm (F'_p)^2.(r'_x+p\ r'_y)+H.r-L.r'_p\equiv c\ F.$$
En \égalant les coefficients des polyn\^omes de degr\é $\leq 3d-4$ en $p$ dans $(*_0)$, on obtient un syst\ème de   $3d-3$ \équations  \à $3d-5$ inconnues qui sont les $d-2$ coefficients $r_j$ de $r$ et les $2d-3 $ coefficients $c_u$ de $c$.

Les $2d-3$ derni\ères \équations\footnote{ordonn\ées dans l'ordre croissant du degr\é auquel correspondent ces coefficients.} forment un syst\ème cram\érien par rapport \à ces inconnues parasites que sont les coefficients    $c_u$, ce qui va nous permettre d'\éliminer ces derniers en r\ésolvant ce syst\ème cram\érien et en reportant la solution dans les $d$ premi\ères \équations qui ne feront plus intervenir alors que les $d-2$ inconnues $r_j$ ainsi que leurs d\ériv\ées. Pour cela,  il va \^etre commode d'utiliser des notations matricielles.

\subsection{ Notations matricielles :}

 Tout d'abord,  aucune confusion  n'\étant \à craindre, on notera souvent de la m\^eme fa\ç\  on un  polyn\^ome 
 $P=\sum_{i=0}^k g_i(x,y)\ p^i$ de degr\é $\leq k$ et le $(k+1)$-vecteur colonne qui lui correspond,  tandis que la multiplication  ${\cal O}_h\ [p] \buildrel{\times P}\over\longrightarrow  {\cal O}_{h+k}\ [p]$  par    $P$ 
                           sera repr\ésent\ée     par la   matrice $M(P)$ ci-dessous \à $h+k+1$ lignes et $h+1$ colonnes\footnote{Il n'est pas interdit que $g_k$ soit nul.} :

$$ P=\begin{pmatrix}g_0 \\
                       g_{1} \\
                         \vdots \\
                        g_k  
                         \end{pmatrix} ,\hskip 2cm
M(P)=\begin{pmatrix}g_0&0&0&\cdots&0\\
                       g_{1}&g_0&0&\cdots&0\\
                         \vdots&\vdots&\ddots&&\\
                         g_k&g_{k-1}&&\ddots&0\\ 
                        0& g_k&g_{k-1}&\cdots&g_0\\ 
                         0&0&g_k &\cdots&g_{1}\\ \vdots&\vdots&\vdots&\ddots&\vdots\\
                         0&0&0&0&g_k
                         \end{pmatrix}.
                        $$
Ainsi $r=\sum_{j=0}^{d-3} r_j(x,y) p^j$ sera repr\ésent\é par un   $(d-2)$-vecteur colonne  $\begin{pmatrix}r_{0}\\ \vdots\\ r_{d-2} \\ r_{d-3} \end{pmatrix}$, de m\^eme que  $r'_x$,  $r'_y$, et $r'_p=<N,r>$, $N$ d\ésignant  la matrice  $N=\begin{pmatrix}
0&1&0&0&\cdots&0\\
0&0&2&0&\cdots&0\\
0&0&0&3&\cdots&0\\
\vdots&\vdots&\vdots&\vdots&\ddots&\vdots\\
0&0&0&0&\cdots&d-3\\
0&0&0&0&\cdots&0\\
\end{pmatrix}$, 
 de taille $(d-2)\times (d-2)$.

                   \n Consid\èrons en particulier  les matrices  
                      
                    $M(F):{\cal O}_{2d-4}\ [p] \longrightarrow  {\cal O}_{3d-4}\ [p]\hbox{ \hskip 1cm
                      de taille $(3d-3)\times (2d-3)$,}$ 
                        
                        $M\bigl((F'_p)^2\bigr) :{\cal O}_{d-2}\ [p] \longrightarrow  {\cal O}_{3d-4}\ [p]\hbox{ \hskip .5cm de taille $(3d-3)\times (d-1)$}$,

  $M(L) :{\cal O}_{d-3}\ [p] \longrightarrow  {\cal O}_{3d-4}\ [p]$ \hbox{ \hskip 1cm de taille $(3d-3)\times (d-2)$}$
  $),
  
 et  $M(H) :{\cal O}_{d-3}\ [p] \longrightarrow {\cal O}_{3d-4}\ [p]\hbox{ \hskip .7cm de taille $(3d-3)\times (d-2)$}$ 
 \hb  $ \hbox{\hskip 5cm   (la premi\ère ligne de $M(H)$ \étant  form\ée de z\éros)}.$

            \n Posons  :
 $I_0=\begin{pmatrix} 1&0&\cdots&0 \\0&1&\cdots&0\\ \vdots &\vdots&\ddots&\vdots\\0&0&\cdots&1 \\0&0&\cdots&0 \end{pmatrix}$ et
$I^0=\begin{pmatrix} 0&0&\cdots&0\\1&0&\cdots&0 \\0&1&\cdots& 0\\
\vdots &\vdots&\ddots&\vdots\\0&0&0&1  \end{pmatrix}$, matrices  de taille $(d-1)\times (d-2)$.  
 \n La relation   $( *_0)$ s'\écrit alors : 
                         $$( *_0)\hskip 2cm <M\bigl((F'_p)^2\bigr) ),(I_0.r'_x+I^0.r'_y)> +<\bigl( M (H)-M(L).N\bigr),r > =<M (F),c>.$$     Cette     relation   se d\écompose en 
  $$ ( *_0)^+\hskip 1cm <M^+\bigl((F'_p)^2\bigr) ),(I_0.r'_x+I^0.r'_y)> +<\bigl( M^+ (H) -M^+(L).N\bigr),r > =<M^+ (F),c>,$$et
   $$( *_0)^-\hskip 1cm <M^-\bigl((F'_p)^2\bigr) ),(I_0.r'_x+I^0.r'_y)> +<\bigl( M ^-(H) -M^-(L).N\bigr),r > =<M ^-(F),c>,$$
   $M^+$ (resp. $M^-$)  d\ésignant les sous-matrices form\ées avec les $d$ premi\ères lignes (resp. les $2d-3$ derni\ères)  des diff\érentes matrices  ci-dessus \à $3d-3$ lignes.  

\n La matrice $M^-(F)$, qui est triangulaire avec des 1 sur la diagonale, est inversible. On en d\éduit :
$$( *_0)^-\hskip .8cm c=<\bigl(M^-(F)\bigr)^{-1},<M^-\bigl((F'_p)^2\bigr) ),(I_0.r'_x+I^0.r'_y)> +<\bigl( M^- (H) -M^-(L).N\bigr),r > >, $$ que l'on reporte dans 
$(*)^+$, de sorte qu'apr\ès cette \élimination de $c$,  $(*_0)$ devient  :
$$<B,(I_0.r'_x+I^0.r'_y)>+<E,r>=0 \hbox { , o\ù}$$

$B=   M^+\bigl((F'_p)^2\bigr)- M^+(F).\bigl(M^-(F)\bigr)^{-1}. M^-\bigl((F'_p)^2\bigr)    $,  

$E=  \bigl( M^+ (H) -M^+(L).N\bigr)-  M^+(F).\bigl(M^-(F)\bigr)^{-1} \bigl( M^- (H)-M^-(L).N\bigr)   $.

\n On v\érifie alors le 

\n {\bf Lemme 2 :} {\it La matrice $B$, qui s'exprime \à l'aide  des coefficients $a_i$ et ne fait pas intervenir leurs d\ériv\ées,    est de rang maximum $d-1$.} 

%\n{\it D\émonstration :}  ....

  \n Soit  $T$ n'importe quelle     matrice\footnote{Si   $B_i$ d\ésigne  la sous-matrice carr\ée de $B$ obtenue en supprimant la $i$-\ème ligne, et si l'indice $i_0$ a \ét\é choisi de fa\ç on que $B_{i_0}$ soit inversible, on peut prendre pour $T$ la matrice $T_{i_0}$ construite \à partir de $(B_{i_0})^{-1}$ en lui ajoutant une colonne de z\éros \à la ${i_0}$-\ème place. C'est ce qu'on fera dans la programmation qui suit.} inverse \à gauche de $B$ (c'est-\à-dire telle que $T.B$ soit la matrice identit\é $Id_{d-1}$) ;  
on en d\éduit   le
 \n {\bf Th\éor\ème 1 : }{\it L'espace des relations ab\éliennes sur $U$ est isomorphe \à l'espace des familles de fonctions $r=(r_j)_j\   \ \ (0\leq j\leq d-3)$, solutions de l'\équation\footnote{Nous essaierons de distinguer, dans la mesure du possible, le signe $'\equiv\ '$ qui marque une \égalit\é entre fonctions ou germes en un point $m\in U$, et le signe $'=\ '$ qui marque une \égalit\é entre valeurs ou   jets  en   $m$.}  :  }
 $$(*) \hskip 1cm <I_0.r'_x+I^0.r'_y> \equiv\ <M ,r>
\hbox {\hskip 2cm o\ù $M=-T.E$.}$$

 \section{Le cas \él\émentaire $d=3$ :} 
 
 Dans ce cas, le degr\é de $r$ est 0 : $r=r(x,y)$ et  $r'_p\equiv 0$, tandis que $M$ est une  matrice $2\times 1$, c'est-\à-dire un  2-vecteur, dont on notera  $M_1$ et $M_2$ les    composantes. Le syst\ème $(*)$ s'\écrit alors : $$r'_x=M_1.r\ , \hskip 1cm r'_y=M_2.r.$$
 
 \n La connexion sur le fibr\é $R_0$ de rang 1 des 0-jets de relations ab\éliennes    (dont $r$ est une section) s'\écrit $\nabla_x r=r'_x-M_1.r$,  $\nabla_y r=r'_y-M_2.r$, et sa    courbure est donc 
 $$k=(M_2)'_x-(M_1)'_y.$$
 C'est la courbure de Blaschke-Dubourdieu\footnote{Voir par exemple [B].}.
 
  \section{Principe de la programmation du calcul de la courbure dans le cas g\én\éral $d\geq 3$ :} 
  
  De fa\c{c}on g\én\érale, on notera $(...)'_{\alpha,\beta}$ la d\ériv\ée partielle $\frac{\partial^{\alpha+\beta} ...}{\partial x^\alpha.\partial y^\beta}$.  
  
  Notons $R_k$ le fibr\é des relations ab\éliennes formelles \à l'ordre $k$ : 
  $R_0$ est un fibr\é de rang $d-2$, dont l'espace des sections $r=(r_j)_{1\leq j\leq d-3}$ peut-\^etre identifi\é \à ${\cal O}_{d-3}\ [p]$. On r\éservera la notation $r'_{\alpha,\beta}$ aux jets des v\éritables  relations ab\éliennes, et l'on   notera plut\^ot $$u_{\alpha,\beta}=(u_{\alpha,\beta})_j \ \ (\alpha+\beta\leq k,\ 0\leq j\leq d-3)$$les \él\éments constitutifs de $R_k$.
 
Les \él\éments de $R_1$ au dessus d'un \él\ément  $u_{0,0}\in R_0$ sont solutions du syst\ème lin\éaire 
  $$<I_0, u_{1,0}>+<I^0,u_{0,1}> =\ <M ,u_{0,0}>$$
  de rang $d-1$, de $(d-1)$ \équations \à $2(d-2)$ inconnues $(u_{1,0})_j,\ (u_{ 0,1})_j$. 

Plus g\én\éralement, et pour $k\leq d-3$,  les \él\éments de $R_{k+1}$ au dessus d'un \él\ément   $e\in  R_{k}$
  sont solutions d'un  syst\ème lin\éaire $\Sigma_{k+1}(e)$ de rang  $(k+1)(d-1)$ (maximum) :
 Si $e=(u_{\alpha,\beta})_{\alpha+\beta\leq k}$,
   le syst\ème $\Sigma_{k+1}(e)$ est constitu\é   des
   
   $k+1$ \équations $(*_{\alpha,\beta})$ \à valeurs  vectorielles  (ou $(k+1)(d-1)$ \équations \à valeurs scalaires), obtenues en d\érivant  l'identit\é  $(*)$ $\alpha$ fois par rapport \à $x$  et $\beta$ fois par rapport \à $y$,  avec $\alpha+\beta=k$,

  \à $(k+2)(d-2)$ inconnues $(u_{\alpha,\beta})_j $,\ ($j=0,\cdots,d-3$), avec cette fois-ci $ \alpha+\beta=k+1 $.
  
  \n Plus pr\écis\ément,  l'\équation $(*_{\alpha,\beta})$  s'\écrit : 
 $$<I_0, u_{\alpha+1,\beta}>+<I^0,u_{\alpha,\beta+1}> =\ \sum_{\gamma=0}^{\alpha}\sum_{\delta=0}^\beta \begin{pmatrix}\alpha\\ \gamma\end{pmatrix}.\begin{pmatrix}\beta\\ \delta\end{pmatrix}< M'_{\alpha-\gamma,\beta-\delta},u _{ \gamma,\delta }>.  $$
On notera   en abr\ég\é  $<M, u>_{\alpha,\beta}$ le second membre de cette  \équation, et l'on posera :
$$M_{\alpha,\beta}^{\gamma,\delta}:=\begin{pmatrix}\alpha\\ \gamma\end{pmatrix}.\begin{pmatrix}\beta\\ \delta\end{pmatrix} M'_{\alpha-\gamma,\beta-\delta}.$$

\n On notera $P_{k+1}$ la partie homog\ène du syst\ème $\Sigma_{k+1}(e)$, constitu\é par la superposition des expressions $<I_0, u_{\alpha+1,k-\alpha}>+<I^0,u_{\alpha,k+1-\alpha}>$ pour $\alpha$ variant de 0 \à $k+1$.

 \n {\it Notation : }   La plupart  des matrices consid\ér\ées   se d\écomposent en blocs $(d-1)\times (d-2)$ construits 
\n - \à partir des matrices $I^0$, $I_0$, pour $P_k$, 

\n - et \à partir de $M$ et de ses   d\ériv\ées partielles successives  pour les seconds membres des syst\èmes. 

\n Si $(...)$ d\ésigne une telle matrice \à $u(d-1)$ lignes et $v(d-2)$ colonnes, on notera $B_{h,\beta}(...)$ le bloc 

- des lignes 
comprises entre $(\beta-1)(d-1)+1 $ et    $\beta(d-1)$,

- et des colonnes comprises entre $(h-1)(d-2)+1 $ et    $h(d-2)$. 

 \n En particulier, 
$B_{\beta,\beta}(P_{k+1})=I_0$, \ \ $B_{\beta,\beta+1}(P_{k+1})=I^0$, \ \ $B_{ \beta, h}(P_{k+1})=0$ si $h\neq \beta,\beta+1$, et l'on d\émontre ais\ément le

\n {\bf Lemme 3 : }{\it La matrice carr\ée   $P:=P_{d-2}$, de taille $(d-1)(d-2)\times (d-1)(d-2)$,  est  inversible.  }

\n Posons $r^{(h)}=(r'_{\alpha,\beta})_{\alpha+\beta=h}$, et   $u^{(h)}=(u_{\alpha,\beta})_{\alpha+\beta=h}$.

 \n Le syst\ème  $\Sigma_{d-2}(e)$, de la forme 
 $$<P , u^{(d-2)}>=<Q,e>,$$ est donc  cram\érien ; sa partie homog\ène, qui ne contient que des z\éros et des 1,  ne d\épend pas de $e=(u^{(k)})_{0\leq k\leq d-3}$. 
 Il permet de d\éfinir la connexion tautologique de H\énaut sur ${\cal E}:=R_{d-3}$, dont les sections \à d\ériv\ée  covariante  nulle s'identifient aux relations ab\éliennes. 
 
 Soit $s=(u_{\alpha,\beta})_{\alpha+\beta\leq d-3}$ une section de ${\cal E}$ :   chaque $u_{\alpha,\beta}=(u_{\alpha,\beta})_j$,  $(0\leq j\leq d-3)$ est  un $(d-2)$-vecteur (ou plus exactement une fonction 
 $$u_{\alpha,\beta}:(x,y)\mapsto u_{\alpha,\beta}(x,y)$$ sur $U$ \à valeurs dans les $(d-2)$-vecteurs).
 
 \n La connexion tautologique est alors d\éfinie par 
 $$\begin{matrix}\bigl(\nabla _x s\bigr)_{\alpha,\beta} &=&\frac{\partial  u_{\alpha,\beta}}{\partial x}&-&u_{\alpha+1,\beta} &\hbox{\hskip 1cm si } \alpha+\beta\leq d-4,\\  
 &&&&&\\
\bigl(\nabla _x s\bigr)_{\alpha,\beta}  &=&\frac{\partial  u_{\alpha,\beta}}{\partial x}&-&
 <P^{-1},<Q,s>\!>_{\alpha+1,\beta}&\hbox{\hskip 1cm si } \alpha+\beta= d-3.\end{matrix},$$
 
 $$\begin{matrix}\bigl(\nabla _y s\bigr)_{\alpha,\beta} &=&\frac{\partial  u_{\alpha,\beta}}{\partial  y}&-&u_{\alpha ,\beta+1} &\hbox{\hskip 1cm si } \alpha+\beta\leq d-4,\\  
  &&&&&\\
\bigl(\nabla _y s\bigr)_{\alpha,\beta}  &=&\frac{\partial  u_{\alpha,\beta}}{\partial y}&-&
 <P^{-1},<Q,s>\!>_{\alpha ,\beta+1}&\hbox{\hskip 1cm si } \alpha+\beta= d-3.\end{matrix},$$
\indent  o\ù $<Q,s>$ d\ésigne le second membre du syst\ème $\Sigma_{d-2}(e)$.

\n {\bf R\éduction du nombre des inconnues :}

A priori, pour d\éfinir  une section $s$ de ${\cal E}$, il faut se donner $(d-1)(d-2)^2/2$ fonctions $(u_{\alpha,\beta})_j$, avec $0\leq \alpha+\beta\leq d-3$ et $0\leq j\leq d-3$. Nous allons montrer que l'on peut en fait r\éduire ce nombre \à $(d-1)(d-2)/2$. Pour cela, nous     proc\éderons  en deux \étapes :

\n - dans la premi\ère, nous r\éduirons d'abord ce nombre \à $(d-2)^2 $ en montrant  que, pour un ordre de d\érivation $k=\alpha+\beta$ donn\é, les vecteurs $u_{\alpha,\beta}$ sont d\étermin\és par l'ensemble des vecteurs $(u_{0,h})_{0\leq h\leq k}$.

 \n - dans la seconde, nous r\éduirons finalement  ce nombre \à $(d-1)(d-2)/2 $ en montrant  que les composantes   $(u_{0,h})_j$ de $u_{0,h}$ pour $d-2-h\leq j\leq d-3$ ne d\épendent que  des vecteurs $(u_{0,k})_{k \leq h-1}$.

%\pagebreak

 \n {\bf Th\éor\ème 2 : }{\it  
  Une section  $s=(u_{\alpha,\beta})_{\alpha+\beta\leq h}$ de $R_h$ est enti\èrement d\éfinie par ses composantes $ v^{(k)}:=u_{0,k}$ pour $k\leq h$.   Plus pr\écis\ément, $u_{\alpha,\beta}$ s'exprime comme une combinaison lin\éaire
  $$u_{\alpha,\beta}=\sum_{k=0}^{\alpha+\beta} E_{\alpha,\beta}^k \ v^{(k)},$$
  o\ù les matrices $ E_{\alpha,\beta}^k$ ne d\épendent que de $M$ et de ses d\ériv\ées,
  avec la formule de r\écurrence $($d'abord sur $\alpha+\beta$, puis sur $\alpha$ pour $\alpha+\beta$ fix\é$)$ : 
$$E_{\alpha+1,\beta}^h=-(J_0.I^0).E_{\alpha,\beta+1}^h+\sum_{ \gamma=0}^\alpha\sum_{\delta=0}^\beta ( J_0.M_{\alpha,\beta}^{\gamma,\delta}).E_{\gamma,\delta}^h, $$
o\ù $J_0$ d\ésigne la matrice de taille $(d-2)\times (d-1)$ consistant \à ajouter une colonne de z\éros \à droite
de la matrice Identit\é $Id_{d-2}$,  avec la convention $E_{\gamma,\delta}^h =0$ si $h>\gamma+\delta$.
}

   \n {\it D\émonstration : } Puisque $J_0.I_0=I_{d-2}$, il suffit d'appliquer la matrice $J_0$ \à chacun des deux membres de l'\équation $(*)_{\alpha,\beta}$ ci-dessus,  pour obtenir le r\ésultat, .

 \n {\bf Corollaire  : }{\it  En particulier,  $Q$ se d\écompose en blocs $B_{\beta,h}(Q)$ de taille $(d-1)\times (d-2)$, avec}  $$B_{\beta,h}(Q) =\sum_{\gamma=0}^{d-3-\beta}\sum_{\delta=0}^\beta  M_{d-3-\beta,\beta}^{\gamma,\delta}.E_{\gamma,\delta}^h \ ,\hskip 1cm  (0\leq \beta\leq d-2, \ 0\leq h\leq d-3) .$$
  \n {\bf Th\éor\ème 3 : }{\it  Pour tout $h$, $h\leq d-3$, $v^{(h)}=(v^{(h)})_{j=0\cdots d-3}$ est enti\èrement d\éfini par   l'ensemble  $\tilde{v}^{(h)}=(v^{(h)})_{j=0\cdots d-3-h}$ de ses composantes $(u _{ 0,h })_j$ obtenues pour $0\leq j\leq d-3-h$, ainsi que par les $\tilde{v}^{(k)}$ pour $k<h$, avec la formule de r\écurrence
  $$v^{(h)}= (-1)^{h}<K^{h},u_{h,0}>+\sum_{k=0}^{h-1}(-1)^{k} <K^{k-1}.(_0 J),<M,u>_{k,h-1-k}> ,$$
  dans laquelle  $(_0 J)$ d\ésigne la matrice de taille $(d-2)\times (d-1)$ consistant \à ajouter une colonne de z\éros 
\à gauche de la matrice Identit\é $I_{d-2}$, et  $K:=(_0 J).I_0$. }
  \n {\it D\émonstration : } On observe que $(_0 J).I^0=I_{d-2}$,\hb  tandis que 
$<K^n, (w_0,w_1,\cdots ,w_{d-3})>=(w_n,w_{n+1},\cdots,w_{d-3},0,\cdots,0)>
 $ pour tout entier $n\leq d-3$, et $K^{d-2}=0$. Par cons\équent,  l'expression  $(-1)^{h}<K^{h},u_{h,0}>$ de la formule ci-dessus n'intervient pas dans le calcul des composantes $v^{(h)}_j$ pour $j>d-3-h$. 
 
 \n En particulier, appliquant $ (_0 J)$ aux deux membres de l'\équation initiale $(*_0)$, on obtient 
 $$ v^{(1)}=-<K,u_{1,0}>+<(_0 J).M,v^{(0)}>,$$
ce qui prouve le th\éor\ème pour $k=1$. 

Plus g\én\éralement, appliquant $(-1)^{k}K^{d-1-h+k}.(_0 J)$ aux deux membres de l'\équation $(*_{k,h-1-k})$ et sommant terme \à terme les \égalit\és obtenues pour $h$ fix\é en faisant varier $k$ de 0 \à $h$, on obtient finalement la formule annonc\ée.
 
Par contre, on peut fixer librement les inconnues $(u_{0,h})_{j}$
  pour   $0\leq j\leq d-3-h$. Des th\éor\èmes 2 et 3 r\ésulte alors le

  \n {\bf Corollaire : }
  {\it 
  
 \n $(i)$ Pour $k\leq d-2$, $R_k$ est un espace fibr\é vectoriel holomorphe de rang $\sum_{j=0}^k (d-2-j)$. 
  
  \n $(ii)$  La projection $R_{d-2}-\!\!>R_{d-3}$ est un isomorphisme de fibr\és vectoriels.  
     
 \n $(iii)$  En particulier, $R_{d-3}$ est un fibr\é de rang $(1/2)(d-1)(d-2)$.}

\n {\bf Num\érotation des inconnues libres, et trivialisation\footnote{Relativement  \à cette  trivialisation, la seule ligne non nulle de la matrice de courbure est la derni\ère.
 Ceci corrtespond   en fait  au cas particulier $n=2$ du th\éor\ème 4 d\émontr\é en appendice. 
} du fibr\é $R_{d-3}$ :}

  Le couple $(h,j)$ tel que $0\leq h\leq d-3$ et  $0\leq j\leq d-3-h$ va recevoir  le num\éro d'ordre 
 $$a=\frac{h(2d-3-h)}{2}+j+1.$$
 On notera alors $\sigma_a$ (ou $s_j^{(h)}$) la section de $R_{d-3}$ telle que  $(u_{0,h})_j=1$, et $(u_{0,h'})_{j'}=0$ si   $(j',h')\neq (j,h)$,  (avec $ h'\leq d-3$ et $j'\leq d-3-h'$). D'apr\ès le th\éor\ème 3 ci-dessus, les sections $\sigma_a$, pour $a$ variant de $1$ \à $(1/2)(d-1)(d-2)$),  trivialisent le fibr\é  $R_{d-3}$.

  Ainsi, par exemple,   les $(j,h)$ sont rang\ées dans l'ordre 
  
$\sigma_1=s_0^{(0)} \hbox{ \hskip 9cm pour } d=3,$

   $\sigma_1=s_0^{(0)},\sigma_2=s_1^{(0)},\sigma_3=s_0^{(1)} \hbox{ \hskip 6cm pour } d=4,$
    
    $\sigma_1=s_0^{(0)},\sigma_2=s_1^{(0)},\sigma_3=s_2^{(0)},\sigma_4=s_0^{(1)},\sigma_5=s_1^{(1)}, \sigma_6=s_0^{(2)} \hbox{\hskip 1.5cm pour } d=5,$ 
   
$ {\rm etc...}$

\n {\bf Num\érotation des lignes  de $Q$, (index\ées par $(h,j)$,  $0\leq h\leq d-3,\ \ 0\leq j\leq d-3$):} :
 
 \n Le couple $(h,j)$ re\c{c}oit le num\éro d'ordre 
 $$a'=h(d-2)+j+1.$$
 (Si $j\leq d-3-h$, $a'=a+\frac{h(h-1)}{2} $).

\n La forme de connexion    $\omega=(\!(\omega_a^b)\!)$ relative \à la trivialisation pr\éc\édente est donn\ée par la formule:
 $$\omega_a^b=\hbox { composante de }\nabla \sigma_a \hbox { sur } \sigma_b,$$
et la forme de courbure $\Omega=(\!(\Omega_a^b)\!)$ correspondante s'en d\éduit par la formule :
 $$\Omega_a^b=d\omega_a^b+[\omega,\omega]_a^b.$$
 
 \section{Programmation sur Maple 8 :}
 
 Merci aux experts en Maple de  pardonner  aux auteurs leurs maladresses en programmation, et de les leur  signaler  !
 
$>$ restart;
 
$>$ with(LinearAlgebra);

\n {\bf Entr\ée des donn\ées :}

({\it  Les donn\ées entr\ées ci-dessous sont celles du 5-tissu de Bol}). 

{\it Entr\ée de d :}

$>$ d:=5;

$>$ interface(rtablesize=(d)*(d)+3);

{\it Entr\ée de F :}

$>$ apply(a,i,x,y); apply(F,x,y,p);

 {\it   Cas implicite : entr\ée directe de F par ses coefficients $a_i$  (avec\ $a_0=1$) ;\hb 
  \indent   dans le cas explicite, entr\ée des $p_i$   : $F=\prod_i(p-p_i).$}
 
$>$ $p1:=-1; p2:=1 ;p3:=y/(x-1); p4:=y/(x+1); p5:=2*x*y/(x^2+y^2-1);$

$>$  F:=(p-p1)*(p-p2)*(p-p3)*(p-p4)*(p-p5);

$>$  F:=collect$(\%,p)$;

$>$ a(0,x,y):=1;

\n {\it A partir de maintenant, l'\écriture du programme ne d\épend plus des donn\ées introduites
(sauf le choix de l'indice i0 ci-dessous, permettant de calculer un inverse \à gauche de $MB$)}

\n {\bf Premi\ère partie : Calcul de la matrice $M$ de l'\équation initiale : \hb $ <  I_0 ,  r'_x> + <I^0 ,  r'_y> =  <M , r> $ :}

$>$ $F_p$:=diff(F,p);

$>$ $F_x$:=collect(diff(F,x),p);

$>$ $F_y$:=collect(diff(F,y),p);

$>$ H:=simplify(diff$(F_p,p)*(F_x+p*F_y)-F_p*diff(F_x+p*F_y,p))$;

$>$ H:=collect(H,p);

$>$ L:=simplify($F_p*(F_x+p*F_y))$;

$>$ L:=collect(L,p);

$>$ for j to d-2 do for i to j-1 do h(i,j):=0 od od;

$>$ for j to d-2 do for i from j to  2*d-1+j do h(i,j):=simplify(coeff(L,p,i-j)) od od;

$>$ for j to d-2 do for i from  2*d+j to 3*d-3 do h(i,j):=0 od od;
 
$>$ ML:=Matrix(3*d-3,d-2,h);

$>$ for j to d-2 do for i to j-1 do k(i,j):=0 od od;

$>$ for j to d-2 do for i from j to  2*d-1+j do k(i,j):=simplify(coeff(H,p,i-j)) od od;

$>$ for j to d-2 do for i from  2*d+j to 3*d-3 do k(i,j):=0 od od;

$>$ MH:=Matrix(3*d-3,d-2,k);

$>$ apply(f,i,j); for j to (2*d-3) do for i  to j-1 do f(i,j):=0 od od;

$>$ for j to (2*d-3) do for i from j to j+d do f(i,j):=coeff(F,p,i-j)  od od;

$>$  for j to (2*d-3) do for i from  j+d+1 to 3*d-3 do f(i,j):=0 od od;

$>$ MF:=Matrix(3*d-3,2*d-3,f);

$>$ Fsup:=DeleteRow(MF,d+1..3*d-3);

$>$ Finf:=DeleteRow(MF,1..d);

$>$ IFinf:=simplify(MatrixInverse(Finf));

$>$ expand($(F_p)^2)$;

$>$ collect  (%,p) $ ;

$>$ apply(g2,i,j); for j to (d-1) do for i  to j-1 do g2(i,j):=0 od od;

$>$ for j to (d-1) do for i from j to (j+2*d-1) do g2(i,j):=coeff$((F_p)^2,p,i-j)$  od od;

$>$  for j to (d-1) do for i from  (j+2*d) to 3*d-3 do g2(i,j):=0 od od;

$>$ MF2p:=simplify(Matrix(3*d-3 ,d-1,g2));

$>$ MB:=simplify(simplify(simplify(DeleteRow(MF2p,d+1..3*d-3))-Fsup.IFinf.simplify(DeleteRow(MF2p,1..d))));

\n {\it Il existe i0 tel que Determinant(DeleteRow(MB,i0)) ne soit pas nul ;
$($modifier i0 si besoin est, en fonction du tissu introduit au d\épart$)$.}
 
$>$ i0:=d-2;

$>$ Rank(DeleteRow(MB,i0));

$>$ IMB1:=factor(factor(MatrixInverse(DeleteRow(MB,i0))));

$>$ apply(t,i,j); for i to d-1 do t(i,i0):=0 od ;for j  to i0-1 do for i to d-1 do t(i,j):=IMB1[i,j] od od; for j from i0+1 to d do for i to d-1 do t(i,j):=IMB1[i,j-1] od od ;

$>$ T:=simplify(simplify(Matrix(d-1,d,t)));

$>$ apply(n,i,j);for i to d-2 do for j to i do n(i,j):=0 od od; for i to d-2 do n(i,i+1):=i od; for i to d-2 do for j from i+2 to d-2 do n(i,j):=0 od od;

$>$ YY:=Matrix(d-2,d-2,n);

$>$ MHL:=simplify(simplify(MH-ML.YY));
 
$>$ MC:=simplify(simplify(simplify(DeleteRow((MHL),d+1..3*d-3))-Fsup.IFinf.simplify(DeleteRow((MHL),1..d))));

$>$ M:=-simplify(T.MC);

$>$ apply(MM,i,j);

$>$ for i to d-1 do for j to d-2 do MM(i,j):=simplify(simplify(M[i,j])) od od;

\n {\bf Deuxi\ème partie : calcul des prolongements de l'\équation initiale :}

{\it Calcul de $M'_{a, b}$ :}

$>$ ds:=proc(u,a,b) description "donne la d\ériv\ée d'ordre sup\érieur"; 

 if (a=0 and b=0) then  u  else
 
   if  (evalf(a)$>$0 and b=0) then simplify(diff(u,x\$a)) else 
   
        if (a=0 and evalf(b)$>$0) then simplify(diff(u,y\$b))  else 
        
             simplify(diff(diff(u,y\$b), \$a))
             
 fi fi fi end proc;

$>$ dM:=proc (a,b) m(a,b):=(i,j)-$>$ds(MM(i,j),a,b); Matrix(d-1,d-2,m(a,b)) end proc;

$>$ apply(vJo,x,y);

for i to d-2 do for j to i-1 do vJo(i,j):=0 od od;

for i to d-2 do   vJo(i,i):=1 od;

for i to d-2 do for j from i+1 to d-1 do vJo(i,j):=0 od od;

$>$ Jo:=Matrix(d-2,d-1,vJo); 

$>$ apply(usup0,i,j);

for j to d-2 do usup0(1,j):=0 od; 

for i from 2 to d-1 do for j to i-2 do usup0(i,j):=0 od od; 

for i from 2 to d-1 do usup0(i,i-1):=1  od;

for i from 2 to d-1 do for j from i to d-2 do usup0(i,j):=0 od od;
 
$>$ apply(uinf0,i,j);

 for j to d-2 do uinf0(d-1,j):=0 od;
 
 for i  to d-2 do 

for j to i-1 do
 
uinf0(i,j):=0 od od;
 
for i  to d-2 do 

uinf0(i,i):=1 od ;

for i  to d-2 do 

 for j from i+1 to d-2 do

 uinf0(i,j):=0 od od;
  
$>$ Isup0:=Matrix(d-1,d-2,usup0);

$>$ Iinf0:=Matrix(d-1,d-2,uinf0);

\n {\it Calcul des matrices E(h,a,b) telles que $r'_{a,b} = \sum_{h=0}^{a+b} < E(h,a,b) , r'_{0,h} > $:}

$>$ apply(E,h,a,b);

 for h from 0 to d-2 do  for n from 0 to h-1 do 
  
  E(n,0,h):=Matrix(d-2,d-2,0)  od od ;
  
  for a from 0 to d-2 do for b from 0 to d-2 do  for h from a+b+1 to d-2 do
  
  E(h,a,b):=Matrix(d-2,d-2,0) od od od ;
  
 for h from 0 to d-2 do E(h,0,h):=IdentityMatrix(d-2) od; for k  from 0 to d-3 do  for a from 0 to k do for h from 0 to k+1 do
 
  E(h,a+1,k-a):= simplify(-Jo.Isup0.E(h,a,k-a+1)+
  
 sum('binomial(a,c)*sum('binomial(k-a,e)*Jo.dM(a-c,k-a-e).E(h,c,e)','e'=0..k-a)','c'=0..a)): od od od;
 
 $>$ for k  from 0 to d-3 do for a from 0 to k do for h from 0 to k+1 do

 print('E'(h,a+1,k-a)=E(h,a+1,k-a)) od od od;
  
\n {\it Calcul des matrices G (h,a,b) telles que  $ <M , r> '_{a,b} = \sum_{h=0}^{a+b} < G (h,a,b) , r'_{0,h}>$ :}

$>$ apply(G,h,a,b);

$>$ for k  from 0 to d-3 do for a from 0 to k do for h from 0 to k do

 G(h,a,k-a):= simplify(simplify(
 
 sum('binomial(a,c)*sum('binomial(k-a,e)*dM(a-c,k-a-e).E(h,c,e)','e'=0..k-a)','c'=0..a))); od od od;

$>$ for k  from 0 to d-3 do for a from 0 to k do for h from 0 to k do

 print('G'(h,a,k-a)=G(h,a,k-a)) od od od;

 \n {\it Expression du syst\ème  $ < P ,  ( r' _{a,b})_{a+b=d-2} >  = <  Q ,  ( r' _{a,b})_{a+b<d-2}  >$ :}

 \n {\it Calcul de la matrice inversible P, partie homog\ène du syst\ème, et matrice inverse IP :}

$>$ apply(u,i,j);

$>$ for a from 1 to d-2 do for i from (a-1)*(d-1)+1 to a*(d-1) do for b from 1 to a-1 do for j from (b-1)*(d-2)+1 to b*(d-2) do 

u(i,j):=0 ; od od od od;

$>$ for a from 1 to d-2 do for i from (a-1)*(d-1)+1 to a*(d-1) do for j from (a-1)*(d-2)+1 to a*(d-2) do 

 u(i,j):=uinf0(i-(a-1)*(d-1),j-(a-1)*(d-2)) od  od od;
 
$>$ for a from 1 to d-2 do for i from (a-1)*(d-1)+1 to a*(d-1) do for j from a*(d-2)+1 to (a+1)*(d-2) do 

u(i,j):=usup0(i-(a-1)*(d-1),j-a*(d-2)) od od od;

$>$ for a from 1 to d-2 do for i from (a-1)*(d-1)+1 to a*(d-1) do  for b from a+2 to d-1 do for j from (b-1)*(d-2)+1 to b*(d-2) do 

 u(i,j):=0 od od od od;
 
$>$ P:=Matrix((d-1)*(d-2),(d-2)*(d-1),u);

$>$ IP:=MatrixInverse(P);

\n {\it Calcul de Q :}

$>$ apply(q,i,j);

$>$ for b from 0 to d-3 do for h from 0 to d-3 do for j from h*(d-2)+1 to (h+1)*(d-2) do for i from b*(d-1)+1 to (b+1)*(d-1) do

 q(i,j):=simplify(G(h,d-3-b,b))[i-b*(d-1),j-h*(d-2)] od od od od ;
 
$>$ Q:=Matrix($(d-2)*(d-1),(d-2)^2,q)$;

$>$ U:=IP.Q;

\n {\it Expression des inconnues  li\ées v(h,j) pour  j variant de d-2-h \à d-3,\hb  en fonction des inconnues libres v(h,j) pour  j variant de 0 \à d-3-h :}

$>$ apply(voJ,x,y);

$>$ for i to d-2 do for j to i do voJ(i,j):=0 od od;

$>$ for i to d-2 do   voJ(i,i+1):=1 od;

$>$ for i to d-2 do for j from i+2 to d-1 do voJ(i,j):=0 od od;

$>$ oJ:=Matrix(d-2,d-1,voJ);

$>$ J:=(oJ).Iinf0;

\n {\it Num\érotation des inconnues libres : $(h,j)\to  a = h*(2*d-3-h) / 2 +j+1 $,   $(0  \leq j \leq d-3-h)$ :}

$>$ apply(hh,a);

$>$ for a to (d-1)*(d-2)/2 do  for h from 0 to d-3 do 

 if h*(2*d-3-h)/2$<$a  and a$<=$(h+1)*(2*d-4-h)/2 then hh(a):=h  fi od od;
 
$>$  for a to (d-1)*(d-2)/2 do print('hh'(a)=hh(a)) od;

$>$ apply(jj,a); for a to (d-1)*(d-2)/2 do jj(a):=a-1-hh(a)*(2*d-3-hh(a))/2 od;

\n {\it Num\érotation de toutes les inconnues  :$ (h,j)\to  a =  h*(d-2) +j+1,\    (0 \leq j  \leq d-3)$ :}

$>$ apply(hhh,a);

$>$ for a to (d-2)*(d-1) do  for h from 0 to d-2 do 

 if h*(d-2)$<$a  and a$<$=(h+1)*(d-2) then hhh(a):=h  fi od od;
 
$>$  for a to (d-2)*(d-1) do print('hhh'(a)=hhh(a)) od;

$>$ apply(jjj,a); for a to (d-2)*(d-1) do jjj(a):=a-1-hhh(a)*(d-2) od;

\n {\it Trivialisation du fibr\é porteur de la connexion par les sections $s(a) =(s(a,h,j))_{h,j}, \ (j\leq d-3-h) $:}

$>$ apply(s,a,h,j);

$>$ for a to (d-1)*(d-2)/2 do for h from 0 to d-3 do  for j from 0 to d-3-h do

 s(a,hh(a),jj(a)):=1 od od od ;

$>$ for a to (d-1)*(d-2)/2 do for h from 0 to d-3 do  for j from 0 to d-3-h do

 if h$<>$hh(a) or j$<>$jj(a) then s(a,h,j):=0 fi od od od ;

\n {\it Expression des inconnues li\ées $(j> d-3-h )$ en fonction des inconnues libres $(j $ au plus $d-3-h)$ :}

$>$ for a to (d-1)*(d-2)/2 do for h from 1 to d-3 do  for j from  d-2-h to d-3 do

s(a,h,j):= sum(' sum(' sum('

 $(-1)^k*(J^k.oJ.G(n,k,h-1-k))[j+1,r]*s(a,n,r-1)',
 'r'=1..d-2) ','k'=0..h-1) ','n'=0..h-1)$  od od od ;
 
$>$ for a to (d-1)*(d-2)/2 do for h from 0 to d-3 do  for j from 0 to d-3 do print('s'(a,h,j)=s(a,h,j)) od od  od; 

$>$ apply(WW,a,h);

$>$ for a to (d-1)*(d-2)/2 do for h from 0 to d-3 do WW(a,h):=Vector(d-2) od od; 

for a to (d-1)*(d-2)/2 do for h from 0 to d-3 do for i to d-2 do 

 WW(a,h)[i]:=s(a,h,i-1) od od od;
  
$>$ for a to (d-1)*(d-2)/2 do for h from 0 to d-3 do 

 print('WW'(a,h)=WW(a,h)) od od;

$>$ apply(WWW,a);

$>$ for a to (d-1)*(d-2)/2 do  WWW(a):=Vector$((d-2)^2)$ od:

$>$ for a to (d-1)*(d-2)/2 do  for i to $(d-2)^2$ do 

WWW(a)[i]:=WW(a,hhh(i))[i-hhh(i)*(d-2)] od  od:

$>$ for a to (d-1)*(d-2)/2 do 

print('WWW'(a)=WWW(a)) od;
 
\n {\bf Troisi\ème partie : expression  de la connexion et calcul de sa courbure :}

({\it Relativement \à la trivialisation pr\éc\édente, la forme de connexion est \écrite $ (Ax)\ dx+(Ay)\ dy$ })

\n {\it D\érivation covariante en x des sections s(a) ; calcul de Ax :}

$>$ apply(Nablax,a,h,j);

$>$  for a to (d-1)*(d-2)/2 do  for h from 0 to d-4 do for j from 0 to d-3 do

  Nablax(a,h,j):=simplify(diff(s(a,h,j),x)-sum('(E(k,1,h).WW(a,k))[j+1]','k'=0..h+1)) od od od;
  
$>$ for a to (d-1)*(d-2)/2 do for h from 0 to d-4 do  for j from 0 to d-3 do

 print('Nablax'(a,h,j)=Nablax(a,h,j)) od od od;

$>$  for a to (d-1)*(d-2)/2 do for j from 0 to d-3 do 

 Nablax(a,d-3,j):=diff(s(a,d-3,j),x)- (U.WWW(a))[(d-3)*(d-2)+j+1] od od;

$>$ for a to (d-1)*(d-2)/2 do for j from 0 to d-3 do

 print('Nablax'(a,d-3,j)=Nablax(a,d-3,j)) od od;
 
$>$ apply(NNx,b,a);

$>$ for b to (d-1)*(d-2)/2 do for a to (d-1)*(d-2)/2 do 

NNx(b,a):=simplify(Nablax(a,hh(b),jj(b))) od od;

$>$ Ax:=Matrix((d-1)*(d-2)/2,(d-1)*(d-2)/2,NNx);

\n {\it D\érivation covariante en y ; calcul de Ay :}

$>$ apply(Nablay,a,h,j);

$>$  for a to (d-1)*(d-2)/2 do for h from 0 to d-4 do  for j from 0 to d-3 do

 Nablay(a,h,j):=simplify(diff(s(a,h,j),y)-s(a,h+1,j)) od od od;

$>$ for a to (d-1)*(d-2)/2 do for h from 0 to d-4 do for j from 0 to d-3 do

 print('Nablay'(a,h,j)=Nablay(a,h,j)) od od od;

$>$  for a to (d-1)*(d-2)/2 do for j from 0 to d-3 do 

  Nablay(a,d-3,j):=simplify$(diff(s(a,d-3,j),y)- (U.WWW(a))[(d-2)^2+j+1])$ od od;
  
$>$ for a to (d-1)*(d-2)/2 do for j from 0 to d-3 do

 print('Nablay'(a,d-3,j)=Nablay(a,d-3,j)) od od;
 
$>$ apply(NNy,b,a);

$>$ for b to (d-1)*(d-2)/2 do for a to (d-1)*(d-2)/2 do 

NNy(b,a):=Nablay(a,hh(b),jj(b)) od od;

$>$ Ay:=Matrix((d-1)*(d-2)/2,(d-1)*(d-2)/2,NNy);

\n {\it Calcul de la courbure :}

$>$ apply(axy,b,a);

for a to (d-1)*(d-2)/2 do for b to (d-1)*(d-2)/2 do 

axy(b,a):=simplify(diff(Nablay(a,hh(b),jj(b)),x)) od od;

$>$ apply(ayx,b,a); for a to (d-1)*(d-2)/2 do for b to (d-1)*(d-2)/2 do 

ayx(b,a):=simplify(diff(Nablax(a,hh(b),jj(b)),y)) od od;
 
$>$ Axy:=simplify(Matrix((d-1)*(d-2)/2,(d-1)*(d-2)/2,axy));

$>$ Ayx:=simplify(Matrix((d-1)*(d-2)/2,(d-1)*(d-2)/2,ayx));
 
$>$ KK:=simplify(simplify(simplify(Axy)-simplify(Ayx))+

simplify(simplify(Ax).simplify(Ay))-simplify(simplify(Ay).simplify(Ax))):
 
\n {\it Dans le cas d'une d\éformation du tissu à l'aide du param\ètre z (sinon KO=KK) :}

$>$ ko:=(i,j)-$>$taylor(KK[i,j],z,1):

$>$ KO:=Matrix((d-1)*(d-2)/2,(d-1)*(d-2)/2,ko);

 \section{Appendice : Concentration de la matrice de courbure  relative \à une trivialisation adapt\ée }
 Nous nous placerons plus g\én\éralement dans le contexte d'un $d$-tissu  calibr\é ordinaire\footnote{On se r\éf\ère \à [DL] pour la terminologie et les notations.} de codimension un dans une vari\ét\é holomorphe de dimension $n$ ($n\geq 2)$. 
 
 Soit $k_0$ l'entier $\geq 2$ tel que  
 $d=(n-1+k_0)!/(n-1)!(k_0)!$ (si $n=2$, $k_0=d-1$).
 
 On d\éfinit alors  une filtration d\écroissante de ${\cal E} =R_{k_0-2}$ en posant : 
 $$F_0({\cal E})={\cal E}\hbox{  et }F_h({\cal E})={\rm Ker}(R_{k_0-2}\to R_{h-1}) \hbox{  pour } 1\leq h\leq k_0-1,$$
 chaque $F_h({\cal E})$ \étant alors un sous-fibr\é vectoriel de rang $\sum_{k=h+1}^{k_0-1}\bigl(d-c(n,h)\bigr)$, 
 o\ù l'on a pos\é 
 $$c(n,h):=(n-1+h)!/(n-1)!h!.$$
 Une trivialisation holomorphe $(\sigma_a)$ de ${\cal E}$ $\Bigl(1\leq a\leq \sum_{h=1}^{k_0-1}\bigl(d-c(n,h)\bigr)\Bigr)$ sera dite \emph{adapt\ée} si, pour tout $h$ $(1\leq h\leq k_0-2)$,  les $(d-c(n,h+1)$ sections $\sigma_a$ telles que 
 $$\sum_{k=1}^{h}\bigl(d-c(n,k)\bigr)< a\leq \sum_{k=1}^{h+1}\bigl(d-c(n,k)\bigr)$$
  engendrent   un suppl\émentaire    $S_{h}$ (n\éc\éssairement holomorphe) de $F_{h+1}({\cal E})$ dans $F_{h}({\cal E})$.
 
 \n {\bf Th\éor\ème 4 :}
{\it  Relativement  \à une   trivialisation ``adapt\ée'',  la matrice de courbure de la connexion tautologique d'un $d$-tissu holomorphe ordinaire calibr\é de codimension un dans une vari\ét\é holomorphe de dimension $n$ $(n\geq 2)$ est concentr\ée\footnote{Ce résultat  probablement d\éj\à bien connu pour $n=2$.} dans les $(n-2+k_0)!/(n-2)!k_0!$ derni\ères lignes $($la derni\ère ligne si $n=2)$, o\ù $k_0$ d\ésigne l'entier $\geq 2$ tel que $d=(n-1+k_0)!/(n-1)!k_0!$.}
\n Cela veut dire que  les lignes qui pr\éc\èdent n'ont que des 0.
 
 \n {\it D\émonstration :} Soit $L=(\alpha_1,\cdots,\alpha_n)$ un multi-indice de d\érivation,  et $L+1_i$ le multi-indice obtenu en augmentant $\alpha_i$ de 1. Supposons la  trivialisation $(\sigma_a)$  adapt\ée. Une    section $s$ de ${\cal E}$ est alors repr\ésent\ée par une famille $(s^{(h)})_h$ d'\él\éments $s^{(h)}\in S_h$, ($s^{(h)}=(s_L)_{|L|=h}$). 
 
 De m\^eme, soit  $\hat s^{(k_0-1)}=(\hat s_L)_{|L|=k_0 -1}$ la famille  des d\ériv\ées partielles   d'ordre $k_0 -1$  des composantes de la projection sur $R_0$  d'une section $\hat s$ de  $R_{k_0-1}$   (composantes relatives \à la trivialisation locale de $R_0$  d\éj\à utilis\ée  pour les d\ériv\ées d'ordre $\leq k_0-2$). 
  
 D'autre part, le tissu \étant ordinaire et calibr\é, 
la projection   $P:R_{k_0-1}\buildrel\cong\over\rightarrow R_{k_0-2}(={\cal E})$ est un isomorphisme de fibr\és vectoriels. Le syst\ème lin\éaire exprimant qu'une section  $\hat s$ de $R_{k_0-1}$ est un $(k_0-1)$-jet de relation ab\élienne se projetant sur une section $s$ de 
$R_{k_0-2}(={\cal E})$ s'\écrit : $$<P,\hat s^{(k_0-1)}>=<Q,s>$$ pour un certain op\érateur lin\éaire $Q$. 
La d\érivation covariante de la connexion tautologique sur ${\cal E}$ est donc d\éfinie par les formules :
 $$(\nabla_i s)_L=\partial_i (s_L)-s_{L+1_i} \hbox{ pour }|L|\leq k_0-3, $$
 $$\hbox{ et }(\nabla_i s)_L=\partial_i (s_L)-<(P^{-1}.Q),s>_{L+1_i} \hbox{ pour }|L|= k_0-2.  $$
 On en d\éduit : 
 $$(\nabla_i\nabla_j s)_L=\partial_i \partial_j (s_L)-\partial_i(s_{L+1_j})-\partial_j( s_{L+1_i})+s_{L+1_i+1_j}  \hbox{ pour }|L|\leq k_0-4,$$
$$\hbox{ et }(\nabla_i\nabla_j s)_L=\partial_i \partial_j (s_L)-\partial_i(s_{L+1_j})-\partial_j( s_{L+1_i})+<(P^{-1}.Q),s>_{L+1_i+1_j}  \hbox{ pour }|L|= k_0-3.$$
Dans les deux cas, l'expression est sym\étrique en $i$ et $j$ : la matrice de courbure ne peut donc avoir de composante non nulle que pour $|L|= k_0-2$. Puisque $S_{k_0-2}$ est de rang $c(n,k_0)-c(n,k_0-1)$, soit $c(n-1,k_0)$, le th\éor\ème est d\émontr\é.

Pour $n=2$, $c(n-1,k_0)$ est \égal \à 1 quel que soit $k_0=d-1$. %En particulier, la trace est le seul coefficient du polyn\^ome caract\éristique de la matrice de  courbure pouvant n'\^etre pas nul. 

%\noindent  [H2] A. H\énaut, Syst\èmes diff\érentiels, nombre de Castelnlication uovo, et rang des tissus de {\bb C}$^n$, Publ. RIMS, Kyoto University, 31(4), 1995, 703-720.

   \n Jean-Paul Dufour, ancien professeur \à l'Universit\é de Montpellier II, \hb
  1 rue du Portalet, 34820 Teyran, France\hb  email : dufourh@netcourrier.com,

  \n Daniel Lehmann, ancien professeur \à l'Universit\é de Montpellier II,\hb   4 rue Becagrun,  30980 Saint Dionisy, France\hb  email : lehm.dan@gmail.com,

\end{document}